\documentclass[oneside,english]{siamltex}
\usepackage[T1]{fontenc}
\usepackage[latin9]{inputenc}
\usepackage{geometry}
\usepackage{array}
\usepackage{multirow}
\usepackage{amsbsy}
\usepackage{amssymb}
\usepackage{graphicx}
\usepackage{esint}

\makeatletter

\providecommand{\tabularnewline}{\\}

\newcommand\eqref[1]{(\ref{#1})}

\usepackage{tikz}
\usetikzlibrary{automata, positioning, arrows}

\@ifundefined{showcaptionsetup}{}{%
 \PassOptionsToPackage{caption=false}{subfig}}
\usepackage{subfig}
\makeatother

\usepackage{babel}
\begin{document}

\title{Variable Viscosity and Density Biofilm Simulations using an Immersed
Boundary Method, Part II: Experimental Validation and the Heterogeneous
Rheology-IBM }

\author{Jay A. Stotsky\thanks{Department of Applied Mathematics, University of Colorado, Boulder, CO 80309-0526} \and Jason
F.~Hammond\thanks{High Power Microwave Division, AFRL, Kirtland AFB, Albuquerque, NM 87116} \and Leonid
Pavlovsky$^{\ddagger}$ \hspace{-3mm} \and Elizabeth J.~Stewart\thanks{Department of Chemical Engineering, University of Michigan, Ann Arbor, MI 48109} \and John
G.~Younger\thanks{Department of Emergency Medicine, University of Michigan Ann Arbor, MI 48109} \and Michael
J.~Solomon$^{\ddagger}$ \hspace{-3mm} \and David M.~Bortz$^*$\thanks{Corresponding author (dmbortz@colorado.edu)}}
\maketitle
\begin{abstract}
The goal of this work is to develop a numerical simulation that accurately
captures the biomechanical response of bacterial biofilms and their
associated extracellular matrix (ECM). In this, the second of a two-part
effort, the primary focus is on formally presenting the \emph{heterogeneous
rheology Immersed Boundary Method} (hrIBM) and validating our model
against experimental results. With this extension of the Immersed
Bounadry Method (IBM), we use the techniques originally developed
in Part I, (Hammond et al. \cite{hammond2014variable}) to treat the
biofilm as a viscoelastic fluid possessing variable rheological properties
anchored to a set of moving locations (i.e., the bacteria locations).
We validate our modeling approach from Part I by comparing dynamic
moduli and compliance moduli computed from our model to data from
mechanical characterization experiments on \emph{Staphylococcus epidermidis
}biofilms. The experimental setup is described in Pavlovsky et al.
(2013) \cite{pavlovsky2013insitu} in which biofilms are grown and
tested in a parallel plate rheometer. Matlab code used to produce
results in this paper will be available at https://github.com/MathBioCU/BiofilmSim.\end{abstract}

\section{Introduction}

The goal of this work is to develop a numerical simulation method
that accurately captures the biomechanical response of bacterial biofilms
and their associated extracellular matrix (ECM). In this second paper,
we show that the model and simulation method developed in part I \cite{hammond2014variable},
can be used to predict material properties of a biofilm and that the
simulated results mimic experimentally measured results. The underlying
mathematical technique is an adaptation of the Immersed Boundary Method
(IBM) that takes into account the finite volume of bacteria, and variable
material parameters found in biofilms whose variation is anchored
to the positions of bacteria in a biofilm. We call this method the
\emph{heterogeneous rheology Immersed Boundary Method} (hrIBM). A
key feature of our results is that the simulations are initialized
with experimentally measured position data providing the locations
of bacteria in live \emph{S. epidermidis} biofilms. This removes ambiguity
about how to represent the biofilm computationally. When using this
data, the bulk physical properties estimated through simulation match
experimental results. We also verify that when using different position
data sets that possess similar spatial statistics, the physical properties
of the biofilm do not change significantly. We also provide quantitative
results on the periodic rotation of suspended aggregates of bacteria
in shear flow. 

In recent years, much work has been done to develop detailed mathematical
models that capture the biomechanical response of bacterial biofilms
to physical changes \cite{alpkvist2007amultidimensional,alpkvist2006threedimensional,danvo2010anexperimentally,hammond2014variable,laspidou2014material}.
In general, the physical properties governing the growth, attachment,
and detachment of a biofilm are dependent on the ECM, a viscous mixture
of polysaccharides and other biological products excreted by bacteria
in the biofilm. The focus of this work is on accurately simulating
the biomechanical response of a biofilm and its associated ECM due
to applied shear stress and shear strain. 

In Section \ref{sec:BiofilmModel}, we provide a brief review of the
classical Immersed Boundary Method (IBM), a well known computational
technique used for the simulation of coupled fluid-structure interactions.
Additionally, we discuss some other IBM based biofilm models, and
explain the adaptations of the IBM that lead to the hrIBM. In Section
\ref{sec:ConvergenceRate}, a description of the numerical properties
and, results from numerical tests showing that the model is convergent
are provided. In Section \ref{sec:Experimental-Validation}, methodologies
for computing relevant material properties from the model are discussed,
and the dynamic moduli and compliance moduli estimated by the model
are compared to experimental data from biofilms grown in a bioreactor.
We observe that these properties do not significantly vary when several
different experimental coordinate data sets with similar spatial statistics
are used. We also compare results of tumbling of bacteria aggregates
suspended in shear flow against theoretical results provided in Blaser
et al. \cite{blaser2002forces}. In Section (\ref{sec:Future-Directions}),
we discuss future research direction and limitations. 

The ability to calculate bulk material properties of a biofilm while
directly incorporating the microscale rheology and connectivity of
the biofilm is the primary contribution of this article. This development
shows that IBM-based models which connect fine scale features such
as models that describe viscoelastic connections between bacteria,
to fluid dynamical models can provide physically accurate results.
From our results, we see that the hrIBM model accurately captures
the elastic component of the biomechanical response of biofilms to
applied stress and strain, and matches experimental trends observed
in the viscous response.

To our knowledge, this work is the first to use a model that accounts
for both the heterogeneous rheological properties and the inter-bacterial
connectivity to compute material properties of a biofilm. Code used
to produce the results obtained in this paper will be available at
https://github.com/MathBioCU/BiofilmSim.

\section{The Biofilm Model\label{sec:BiofilmModel}}

In this section, we discuss some previous biofilm models and explain
the alterations of the classical IBM that lead to the hrIBM. In Section
\ref{sub:The-Biofilm-Model}, we introduce our biofilm model. In our
model, we couple a spring model of the inter-bacteria links in the
biofilm with fluid motion through the biofilm to treat the biofilm
as a multicomponent viscoelastic material. On the level of our simulations,
both the fluid-structure interactions of the bacteria and the surrounding
fluid, and the interconnectedness of bacteria in the biofilm play
a major role.

\subsection{Previous IBM Based Biofilm Models}

In recent years, a number of different approaches to IBM-based biological
material models have been developed. One such biofilm model can be
found in Luo et al. \cite{luo2008animmersedboundary}. In their model,
they couple an immersed viscoelastic structure to the fluid flow in
an immersed boundary type formulation. However, the fluid equations
are solved separately from the equations governing the motion of the
immersed viscoelastic solid and then coupled together at a physical
interface. Our model builds upon this work by eliminating the need
for an explicit interface since biofilms frequently do not have well
defined fluid-structure interfaces. Another approach to capturing
the viscoelastic nature of biofilms with an the immersed boundary
method is through the choice of viscoelastic model used for the links
between bacteria. The choice of model can be used to affect the value
of the external force density, $\mathbf{f}$, in the Navier-Stokes
equations, (\ref{eq:N-S}). This type of strategy was used first by
Bottino in \cite{bottino1998modeling} to model general viscoelastic
connections in actin cytoskeleton of ameboid cells and also by Dillon
and Zhuo in \cite{zhuo2010usingthe} to model sperm motility. 

IBM-based models can be found in Alpkvist and Klapper, \cite{alpkvist2006threedimensional}
and in Dan Vo et al. \cite{danvo2010anexperimentally}. In these models,
an IBM is used directly to couple the forces between connected bacteria
with fluid motion. Additionally, some validation results are performed
to show that properties such as the recovery and relaxation times
of a biofilm can be modeled with such a method. Our model also builds
on these by including spatially variable rheological properties which
are an important structural feature of biofilms. Additionally, we
note the recent work of Sundarson et al. \cite{sudarsan2015simulating}
in which an IBM model is used to model detachment of biofilms.

We would also like to point out that detailed explanations of the
IBM can be found in \cite{peskin1977numerical,peskin2002theimmersed,zhu2003interaction,zhuo2010usingthe},
and additional IBM-based biofilm models can be found in \cite{alpkvist2007amultidimensional,danvo2010anexperimentally,hammond2014variable}.

\subsection{The Biofilm Model\label{sub:The-Biofilm-Model}}

The model we use is comprised of two sets of equations; those that
model fluid flow through the biofilm, and those that model motions
and forces experienced by each bacteria cell in the biofilm. These
equations are listed in (\ref{eq:N-S})-(\ref{eq:visc_transfer}).
\begin{eqnarray}
\rho(\mathbf{x},\, t)\left(\mathbf{u}_{t}+\left(\mathbf{u}\cdot\nabla\right)\mathbf{u}\right)=-\nabla P+\nabla\cdot\mu(\mathbf{x},\, t)\left(\nabla\mathbf{u}+(\nabla\mathbf{u})^{T}\right)+\mathbf{f}(\mathbf{x},\, t)\label{eq:N-S}\\
\nabla\cdot\mathbf{u}=0\label{eq:Continuity_Eq}\\
\mathbf{U}(\mathbf{X}(s,\, t),\, t)=\int_{\Omega}\mathbf{u}(\mathbf{x},\, t)\,\delta(\mathbf{X}(s,\, t)-\mathbf{x})\, d\mathbf{x\;}s=1,2,...N\label{eq:LagrangeVel}\\
\frac{\partial\mathbf{X}(s,\, t)}{\partial t}=\mathbf{U}(\mathbf{X}(s,\, t),\, t)\label{eq:LagrangePosition}\\
\mathbf{F}(\mathbf{X}(s,\, t),\, t)=\mathcal{F}(\mathbf{X}(s,\, t),\mathcal{\, P})\label{eq:LagrangeForce}\\
\mathbf{f}(\mathbf{x},\, t)=\frac{1}{d_{0}^{3}}\int_{\Omega}\mathbf{F}(\mathbf{X}(s,\, t),\, t)\,\hat{\delta}(\mathbf{X}(s,\, t)-\mathbf{x},\,\omega)\, d\mathbf{X}\label{eq:Smoothed_Delta_Transfer}\\
\rho(\mathbf{x},\, t)=\rho_{0}+\min\left\{ \int_{\Omega}\omega^{3}\rho_{b}\,\hat{\delta}(\mathbf{X}(s,\, t)-\mathbf{x},\,\omega)\, d\mathbf{X},\,\rho_{b}\right\} \\
\mu(\mathbf{x},\, t)=\mu_{0}+\min\left\{ \int_{\Omega}\omega^{3}\mu_{b}\,\hat{\delta}(\mathbf{X}(s,\, t)-\mathbf{x},\,\omega)\, d\mathbf{X},\,\mu_{b}\right\} .\label{eq:visc_transfer}
\end{eqnarray}
The same set of equations was also used in Part I, \cite{hammond2014variable}
and are reproduced here for the convenience of the reader. The quantities
appearing in these equations are listed in Table \ref{tab:List-of-terms}. 

\begin{table}
\begin{centering}
\begin{tabular}{|c|>{\raggedright}m{3in}|}
\hline 
Symbol & Definition\tabularnewline
\hline 
\hline 
$\mathbf{u}$, $P$, $\mathbf{F}$ & Eulerian velocity, pressure and force density\tabularnewline
\hline 
$\rho(\mathbf{x},\, t)$, $\mu(\mathbf{x},\, t)$ & spatially (and temporally) varying density and viscosity\tabularnewline
\hline 
$d_{0}^{3}$, $N$ & Average volume associated with a cell and its surroundings, number
of bacteria in the domain\tabularnewline
\hline 
\multirow{1}{*}{$\mathbf{X}(s,\, t)$, $\mathbf{U}(\cdot,\, t)$, $\mathbf{F}(\cdot,\, t)$} & Bacteria position, velocity, and force density, labelled by Lagrangian
coordinate, $s$\tabularnewline
\hline 
$\rho_{0}$, $\mu_{0}$ & Density and viscosity of pure water\tabularnewline
\hline 
$\rho_{b}$, $\mu_{b}$  & Density and viscosity of the biofilm at the center of mass of each
bacteria\tabularnewline
\hline 
$\mathcal{F}(\cdot,\,\mathcal{P})$ & Function that determines the force associated with each bacteria based
on a constitutive viscoelastic model, $\mathcal{P}$\tabularnewline
\hline 
$\delta(\cdot)$ & The Dirac delta function\tabularnewline
\hline 
\multirow{1}{*}{$\hat{\delta}(\cdot,\,\omega)$} & A smoothed approximation of the Dirac delta function. The 2nd argument,
$\omega$ is a hydrodynamic parameter corresponding to the radius
of a bacterium.\tabularnewline
\hline 
\end{tabular}
\par\end{centering}

\protect\caption{\label{tab:List-of-terms}List of terms in governing equations. We
attempt to use standard notations when possible.}
\end{table}

Since individual bacteria are not assumed to have infinitessimal volume
at the scale of our simulations, the Lagrangian quantities; $\mathbf{X}$,
$\mathbf{F}$, and $\mathbf{U}$ correspond to measurements taken
at the center of mass of each bacterium.  As described in Section
\ref{sub:The-Immersed-Boundary-Method}, the second argument, $\omega$,
of the smoothed Dirac $\delta$ function, $\hat{\delta}(\cdot,\,\omega)$
determines a region of support for the smoothed $\delta$ function
. The choice of $\hat{\delta}(\cdot,\,\omega)$ govern how the mass
density, viscosity, and force density vary around each bacterium. 

Additionally, since in each simulation, there is a fixed number, $N$
of bacteria in the computational domain which is independent of the
grid spacing, $h$, we can write the integrals in equations (\ref{eq:Smoothed_Delta_Transfer})-(\ref{eq:visc_transfer})
as summations of the form 
\[
\mathbf{f}(\mathbf{x},\, t)=\frac{1}{d_{0}^{3}}\sum_{s=1}^{N}\mathbf{F}(\mathbf{X}(s,\, t),\, t)\,\hat{\delta}(\mathbf{X}(s,\, t)-\mathbf{x},\,\omega)
\]
\[
\rho(\mathbf{x},\, t)=\rho_{0}+\min\left\{ \sum_{s=1}^{N}\omega^{3}\rho_{b}\,\hat{\delta}(\mathbf{X}(s,\, t)-\mathbf{x},\,\omega),\,\rho_{b}\right\} 
\]
\[
\mu(\mathbf{x},\, t)=\mu_{0}+\min\left\{ \sum_{s=1}^{N}\omega^{3}\mu_{b}\,\hat{\delta}(\mathbf{X}(s,\, t)-\mathbf{x},\,\omega),\,\mu_{b}\right\} .
\]
These summations are slightly different than those used on part I.

By using the IBM as a basis for our biofilm model, we are able to
avoid treating the biofilm as a two phase fluid with a distinct bulk
fluid region and a distinct biofilm region. Instead, the use of variable
rheological properties over the entire domain couples the biofilm
and bulk fluid motions as a single viscoelastic material.

\subsection{The Heterogeneous Rheology Immersed Boundary Method\label{sub:The-Immersed-Boundary-Method}}

We will now describe our reasoning behind equations (\ref{eq:N-S})-(\ref{eq:visc_transfer}).
In our model, we extend the IBM to account for the fixed, finite size
of bacteria and allow for variable physical properties that are anchored
to a moving Lagrangian mesh (i.e. the bacteria positions). We denote
this approach the \emph{heterogeneous rheology IBM} (hrIBM). 

The original IBM was first developed as a means of solving fluid-structure
interaction problems in cardiology and is applicable to problems with
moving, irregularly shaped boundaries \cite{peskin2002theimmersed,peskin1977numerical}.
With the IBM, the fluid velocity fields and pressure are usually solved
for on a fixed, Eulerian grid and the movement of the boundaries due
to fluid motion is tracked by a moving Lagrangian mesh. As material
boundaries are deformed, a constitutive model is used to determine
the force density exerted by the boundary on the fluid around each
Lagrangian point. The Lagrangian force density field is then transferred
to an Eulerian force density field through the use of a discrete approximation
of the following identity,
\begin{equation}
\mathbf{f}(\mathbf{x},\, t)=\int_{\Omega}\mathbf{F}(\mathbf{X}(\mathbf{q},\, t))\,\delta(\mathbf{X}(\mathbf{q},\, t)-\mathbf{x})\, d\mathbf{q}\label{eq:Euler_Lagrange_Transfer}
\end{equation}
where $\delta(\cdot)$ is the Dirac $\delta$ function. In our biofilm
model, the transfer of quantities from the Lagrangian to the Eulerian
grid is done with a smoothed $\delta$ function that differs from
the standard choices used in most IBM literature (see (\ref{sub:The-Biofilm-Model})).
The effect of the Eulerian force term on the velocity and pressure
fields is then found by solving the Navier-Stokes equations, (\ref{eq:N-S})
with appropriate boundary conditions.

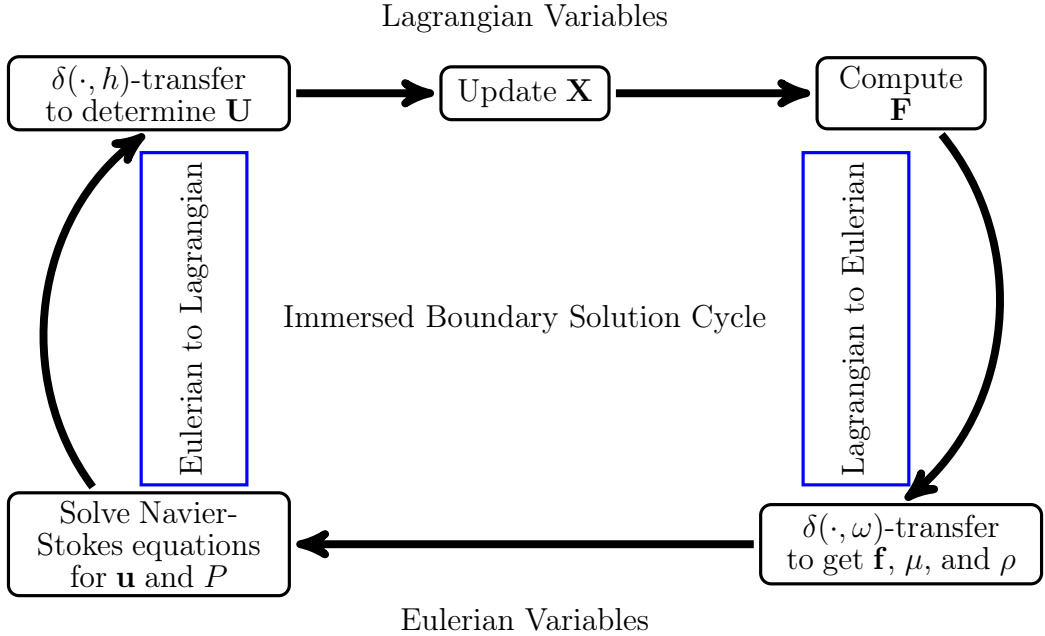
\begin{figure}
\begin{centering}
\thispagestyle{empty}
\tikzset{ %Define standard arrow tip    
	>=stealth',     
%Define style for boxes 
   squ/.style={rectangle, rounded corners, draw=black, very thick, minimum height=2em, text centered}, 
squg/.style={rectangle, draw=blue, very thick, minimum height=4em, rotate=90, text centered}, 
squh/.style={rectangle, draw=blue, very thick, minimum height=4em, rotate=-90, text centered},  
ELt/.style={ rectangle, draw=black, very thick, minimum height=3em, minimum width=10em, fill=blue, opacity=0.0, text centered, text opacity=1},
    % Define arrow style     
pil/.style={->,thick,shorten <=2pt,shorten >=2pt,} }
 \begin{tikzpicture}
%\pgftext { \includegraphics[width=1\textwidth]{biofilmimage.png}} at (0pt, 0pt);
\node[squ] (A) at (0,0)   [rectangle ,text width =3.5cm, draw]{\large Solve Navier-Stokes equations for $\mathbf{u}$ and $P$}; 
\node[squ] (B) at (0,6) [ rectangle,  text width =3.5cm, draw]{\large$\delta(\cdot, h)$-transfer to determine $\mathbf{U}$};
\node[squ](C) at (5,6) [ rectangle, text width=2cm, draw] {\large Update $\mathbf{X}$ };
\node[squ] (D) at (10,6) [rectangle, text width =2cm, draw] {\large Compute $\mathbf{F}$}; 
\node[squ] (E) at(10,0) [rectangle, text width=3.5cm, draw]{\large $\delta(\cdot, \omega)$-transfer to get $\mathbf{f}$, $\mu$, and $\rho$}; \node[squg] (EL) at (0.6,3) {\large Eulerian to Lagrangian}; 
\node[squg] (LE) at (9.4,3){\large Lagrangian to Eulerian};
\node[ELt] (Eul) at (5,-1){\large Eulerian Variables};
\node[ELt] (Lag) at (5,7){\large Lagrangian Variables};
\node[ELt] (Def) at (5,3){\large Immersed Boundary Solution Cycle};
\path  (A)	edge[pil, bend left=45, line width=3pt] (B.south) (B)	edge[pil, bend left=0, line width=3pt] (C.west) (C) edge[pil, bend left=0, line width=3pt] (D.west) (D) edge[pil, bend left=45, line width=3pt] (E.north) (E) edge[pil, bend left=0, line width=3pt] (A.east);
\end{tikzpicture} 
\par\end{centering}

\protect\caption{The coupling between the Eulerian and Lagrangian variables in the
hrIBM is shown here. The Eulerian and Lagrangian variables are coupled
by the computation of $\mathbf{U}$ from $\mathbf{u}$, and the computation
of $\mathbf{f}$, $\mu,$ and $\rho$ from $\mathbf{F}$ and $\mathbf{X}$.
The IBM is a widely applicable method in part because it allows for
a great variety of fluid solvers and solid structural models to be
coupled through $\delta$ function transfer identities.}
\end{figure}

In the original IBM, after discretizing, the integration in (\ref{eq:Euler_Lagrange_Transfer})
is carried out by computing a sum of the form 
\begin{equation}
\mathbf{f}(\mathbf{x}_{i},\, t_{j})=\sum_{k}^{N}\mathbf{F}(\mathbf{X}_{k},\, t_{j})\,\hat{\delta}(\mathbf{X}_{k}-\mathbf{x}_{i},\, h)\, h^{3},\label{eq:EulerLagrangeForce}
\end{equation}
where the Eulerian and Lagrangian forces are evaluated at the Eulerian
and Lagrangian grid points respectively, and $\hat{\delta}(\cdot,\, h)$
is a discrete approximation of the Dirac delta function that has compact
support related to the grid spacing parameter $h$. With the IBM,
the discrete approximation is chosen such that as $h\rightarrow0$,
$\hat{\delta}(\mathbf{r},\, h)\rightarrow\delta(\mathbf{r})$. This
makes sense for fluid structure interactions involving fluid-solid
boundaries that have infinitessimal thickness, and thus zero volume.
In biofilm modeling, each Lagrangian point corresponds to the center
of mass of a bacterium which has finite dimensions. Therefore, we
use a smoothed version of the standard discrete $\delta$ function
that has a fixed region of support, independent of the grid spacing,
which is governed by a radial parameter, $\omega$. 

In our model, we use a smoothed discrete Dirac $\delta$ approximation
of the form,
\begin{equation}
\hat{\delta}(\mathbf{x},\,\omega)=\frac{1}{\omega^{3}}\phi\left(\frac{x}{\omega}\right)\phi\left(\frac{y}{\omega}\right)\phi\left(\frac{z}{\omega}\right)
\end{equation}
with $\phi(r)$ as defined in \cite{peskin2002theimmersed} by 
\[
\phi(r)=\left\{ \begin{array}{cc}
\frac{1}{8}\left(5-2|r|-\sqrt{-7+12|r|-4|r|^{2}}\right) & 1\leq|r|\leq2\\
\frac{1}{8}\left(3-2|r|+\sqrt{1+4|r|-4|r|^{2}}\right) & 0\leq|r|\leq1\\
0 & |r|>2
\end{array}\right.
\]

This is chosen because it most closely satisfies the unity and first-moment
conditions described below for the values of $\omega$ we use. If
$\omega=h$, the standard discrete $\delta$ functions seen in IBM
literature is obtained. For this work, we assume that the bacteria
are spherical and thus $\omega$ is understood as a hydrodynamic radius.
We also note that extensions to this formalism will allow for the
treatment of nonspherical bacteria. Thus, $\omega$ may be though
of more generally as a shape parameter. 

With $\hat{\delta}(\mathbf{x},\, h)$, the \textit{unity condition,}
\begin{equation}
\sum_{\mathbf{x}\in\mathcal{G}_{h}}\hat{\delta}(\mathbf{x}-\mathbf{X},\, h)\thinspace h^{3}=1,\:\forall\mathbf{X},\label{eq:unityConditionDirac}
\end{equation}
and \textit{first-moment condition,}
\begin{equation}
\sum_{\mathbf{x}\in\mathcal{G}_{h}}(\mathbf{x}-\mathbf{X})\,\hat{\delta}(\mathbf{x}-\mathbf{X},\, h)h^{3}=0,\:\forall\mathbf{X},\label{eq:1stMomentConditionDirac}
\end{equation}
are both satisfied. With a grid-independent choice for $\omega,$
these properties are only satisfied approximately. However, we do
see that in the limit as $h\rightarrow0$, greater than $O(h^{2})$
convergence in $\hat{\delta}(\mathbf{r},\,\omega)$ to equations (\ref{eq:unityConditionDirac})
and (\ref{eq:1stMomentConditionDirac}) is observed. 

Highly heterogeneous viscosity and moderately heterogeneous density
are common characteristics of biofilms. Although IB methods with
variable density have existed for some time (see \cite{zhu2003interaction}),
the incorporation of spatially variable viscosity in the IBM is an
area that has yet to be well developed. We do, however note the recent
publications by Fai et al. \cite{fai2013immersed,fai2014immersed}
in which an IBM capable of solving problems with variable viscosity
and density is used to model the motion of red blood cells flowing
in capillaries. When modeling red blood cells, the viscosity exhibits
a ``jump'' discontinuity between the blood plasma and the intracellular,
hemoglobin-containing fluid of a red blood cell. Thus, their model
is designed to capture the dynamics of two interacting fluids with
different rheological properties separated by a deformable membrane.
In our case, there do not exist well defined boundaries and thus,
$\delta(\cdot,\,\omega)$ is adjusted to reflect this. 

In biofilms, the spatial variance of material properties is localized
around the position of each bacterium, while in fluid far away from
any bacteria, the physical properties are those of the bulk fluid.
This localization of the variation in material properties allows the
spatial variation in density and viscosity to be found by using a
smoothed $\delta$-function integration similar to that used to compute
the Eulerian force field. We define an effective viscosity, $\mu_{b}$
and an effective density, $\rho_{b}$ and assume that at the center
of mass of each bacteria, the viscosity and density are $\rho(\mathbf{X}_{i},\, t)=\rho_{b}$,
and $\mu(\mathbf{X}_{i},\, t)=\mu_{b}$. Defining $\mu_{0}$ to be
the viscosity of the bulk fluid, in this case water, the viscosity
at any Eulerian grid point can be calculated as: 
\begin{equation}
\mu(\mathbf{x})=\min\left\{ \mu_{0}+\int_{\Omega}\omega^{3}(\mu_{b}-\mu_{0})\,\hat{\delta}(\mathbf{x}-\mathbf{X}(s),\,\omega)\, d\mathbf{X},\,\,\mu_{b}\right\} .\label{eq:viscosity_function}
\end{equation}
A similar formula exists for the spatial variation of density. A summation
sign is used instead of an integral since the number of bacteria,
$N$, is fixed and independent of the mean Lagrangian mesh spacing
$d_{0}$. 

In this model, we indirectly take into account the fluid volume displacement
caused by the presence of the bacteria. We treat the localized high
viscosity around each bacteria as an effective viscosity that accounts
for both the displaced fluid volume and the increased viscosity near
the bacteria surface \cite{gaboriaud2008coupled}. Extensions based
on changing our choice for $\hat{\delta}(\mathbf{x},\,\omega)$ could
possibly allow for a more precise computation volume displacement
into the model. The bacteria \emph{S. epidermidis} is known to have
a diameter of approximately $0.5-1.0\,\mu m$ \cite{foster1996staphylococcus}
, thus we choose $\omega$ such that the viscosity halo around each
bacterium is a little greater than $1\,\mu m$ in our simulations. 

\begin{figure}
\begin{centering}
\subfloat[]{\begin{centering}
\includegraphics[bb=0bp 0bp 362bp 632bp,width=0.2\paperwidth]{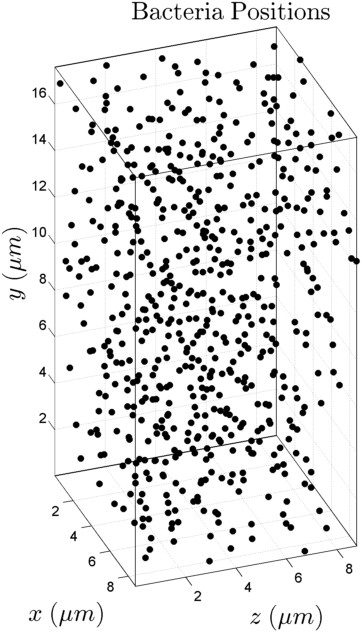}
\par\end{centering}

}\quad{}\subfloat[]{\begin{centering}
\includegraphics[bb=0bp 0bp 360bp 630bp,width=0.2\paperwidth]{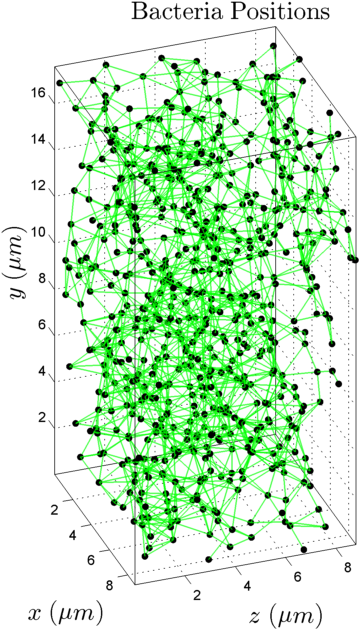}
\par\end{centering}

}\quad{}\subfloat[]{\begin{centering}
\includegraphics[bb=0bp 0bp 362bp 632bp,width=0.2\paperwidth]{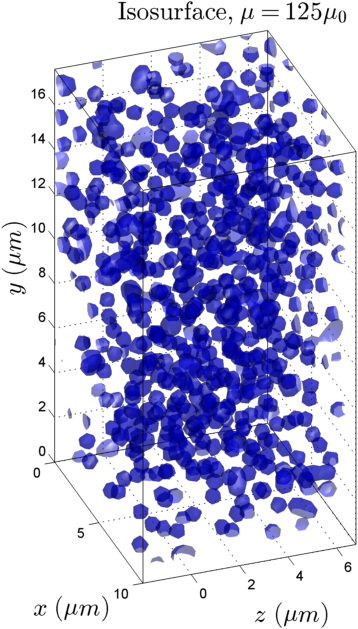}
\par\end{centering}

\centering{}}
\par\end{centering}

\protect\caption{\label{fig:Biofilm_Geom}a) Shows the 3D locations of bacteria from
experimental biofilm data. b) Each line represents a viscoelastic
connection between two bacteria. Bacteria connected if they are within
$1.62\mu m$ of each other. c) A viscosity isosurface of the same
biofilm. The maximum viscosity is $250\mu_{0}$ where $\mu_{0}$ is
the viscosity of water. The isosurface is the surface defined by $\mu(\mathbf{x})=125\mu_{0}.$ }
\end{figure}

\section{Numerical Methods\label{sec:ConvergenceRate}}

The numerical methods we use are based on those originally discussed
in Hammond et al. \cite{hammond2014variable}. We summarize them here
for convenience and also provide convergence results. To approximate
solutions to equations (\ref{eq:N-S})-(\ref{eq:visc_transfer}),
we use a projection method similar to that used in Zhu et al.\cite{zhu2003interaction}.
The solution scheme uses an implicit Euler solver to update an intermediate
velocity profile at each time step and is expected to be $\mathcal{O}(\Delta t)$
convergent. To discretize the domain, we use a uniform finite difference
discretization with equal spacings in the $x$, $y$, and $z$ directions.
The spatial derivatives are approximated with 2nd order, centered
finite differences.

\subsection{Numerical Algorithm}

At each time step, the following quantities must be updated: $\mathbf{u}$,
$\mathbf{U}$, $P$, $\mathbf{F}$, $\mathbf{f}$, $\mathbf{X}$,
$\mu$, and $\rho$. To improve numerical accuracy, we nondimensionalize
the problem with the following choices:

\[
\mathbf{\hat{\mathbf{u}}}=\frac{\mathbf{u}}{u_{0}}\,\,\,\,\hat{P}=\frac{P}{P_{0}}\,\,\,\,\mathbf{\hat{f}}=\frac{\mathbf{f}}{f_{0}}\,\,\,\,\hat{\mu}=\frac{\mu}{\mu_{0}}\,\,\,\,\hat{\rho}=\frac{\rho}{\rho_{0}}\,\,\,\,\mathbf{\hat{\mathbf{x}}}=\frac{\mathbf{x}}{L}\,\,\,\,\hat{t}=\frac{t}{t_{0}}
\]
and also introduce the following nondimensional parameters:
\[
Re=\frac{\rho_{0}Lu_{0}}{\mu_{0}}\,\,\,\, St=\frac{L}{t_{0}u_{0}}\,\,\,\,\, C_{1}=\frac{P_{0}}{\rho_{0}u_{0}^{2}}\,\,\,\,\, C_{2}=\frac{f_{0}L}{\rho_{0}u_{0}^{2}}.
\]

As is standard terminology, $Re$ is the Reynold's number, $St$ the
Strouhal number, and $C_{1}$ and $C_{2}$ are additional constants.
Additionally, we define $d_{0}^{3}$ to be the average Lagrangian
volume element as described in Part I \cite{hammond2014variable}.
For convenience, we will now assume that all quantities are nondimensional
unless otherwise stated. The values of the constants we use are listed
in Table \ref{tab:Values-of-Physical-Constants} and the motivation
for these values is discussed in Part I.

\begin{table}
\begin{centering}
\begin{tabular}{|c|c|}
\hline 
Quantity & Value\tabularnewline
\hline 
\hline 
$P_{0}$ & $1\, Pa$\tabularnewline
\hline 
$\mu_{0}$ & $1\cdot10^{-3}\, Pa\cdot s$\tabularnewline
\hline 
$\rho_{0}$ & $998\, kg/m^{3}$\tabularnewline
\hline 
$L$ & $10^{-5}\, m$\tabularnewline
\hline 
$Re$ & $\mathcal{O}(10^{-3})$\tabularnewline
\hline 
$St$ & $\mathcal{O}(10^{-2})$\tabularnewline
\hline 
$C_{1}$ & $\mathcal{O}(10^{4})$\tabularnewline
\hline 
$C_{2}$ & $\mathcal{O}(1)$\tabularnewline
\hline 
$t_{0}$ & $1\, s$\tabularnewline
\hline 
$f_{0}$ & $1\, N/m^{3}$\tabularnewline
\hline 
$u_{0}$ & $\mathcal{O}(10^{-4})$ (varies)\tabularnewline
\hline 
$d_{0}$ & $0.159\, L$\tabularnewline
\hline 
$\rho_{b}$ & $0.12\rho_{0}$\tabularnewline
\hline 
$\mu_{b}$ & $250\mu_{0}$\tabularnewline
\hline 
$F_{max}$ & $1.3223\cdot10^{-9}$\tabularnewline
\hline 
Connection Distance & $0.162\, L$\tabularnewline
\hline 
Radial Parameter, $\omega$ & $0.033\cdot L$\tabularnewline
\hline 
\end{tabular}
\par\end{centering}

\protect\caption{\label{tab:Values-of-Physical-Constants}Values of Physical Parameters
and Nondimensional constants used in simulations}
\end{table}

As is standard practice in IBM algorithms, we uncouple the Eulerian
variable updates and Lagrangian variable updates for computational
reasons. At each time step, we use a projection-based solver to solve
the Navier-Stokes equation for $\mathbf{u}$ and $P$. We define $\mathbf{G}_{h}$
a discrete gradient operator, and $\mathbf{D}_{h}$ a discrete divergence
operator, and use the following projection method to obtain $\mathbf{u}$
and $P$:
\begin{enumerate}
\item Solve for $\mathbf{u}^{\ast}$ 
\[
\rho^{(n-1)}\left(St\,\frac{\mathbf{u}^{\ast}-\mathbf{u}^{(n-1)}}{\Delta t}+\frac{1}{2}\left(\mathbf{u}^{(n-1)}\cdot\mathbf{D}_{h}(\mathbf{u}^{(n-1)})+\mathbf{D}_{h}\left(\mathbf{u}^{(n-1)}\mathbf{u}^{n-1)}\right)\right)\right)
\]
\[
=\frac{1}{Re}\,\mathbf{D}_{h}\left[\mu^{(n-1)}\left(\mathbf{G}_{h}(\mathbf{u}^{\ast})+(\mathbf{G}_{h}(\mathbf{u}^{\ast}))^{T}\right)\right]+C_{2}\,\mathbf{f}^{(n-1)}
\]

\item Solve for $P^{(n)}$
\[
\mathbf{D}_{h}\left(\frac{1}{\rho^{(n-1)}}\mathbf{G}_{h}P^{(n)}\right)=\left(\frac{St}{C_{1}}\right)\frac{\mathbf{D}_{h}(\mathbf{u}^{\ast})}{\Delta t}
\]

\item Compute $\mathbf{u}^{(n)}$
\[
\mathbf{u}^{(n)}=\mathbf{u}^{\ast}-\left(\frac{C_{1}}{St}\right)\frac{\Delta t}{\rho^{(n-1)}}\mathbf{G}_{h}(P^{(n)})
\]

\end{enumerate}
In steps 1 and 2, full multigrid solvers and multigrid preconditioned
conjugate gradient solvers are used to find $\mathbf{u}^{\ast}$ and
$P^{(n)}$. After obtaining the updated velocity and pressure, the
Lagrangian velocity and position updates follow,
\[
\mathbf{U}^{(n)}=\sum_{h\in\mathcal{G}_{h}}\mathbf{u}^{(n)}\hat{\delta}(\mathbf{x}_{h}-\mathbf{X}^{(n-1)},\, h)\, h^{3}
\]

\[
\mathbf{X}^{(n)}=\mathbf{X}^{(n-1)}+\frac{\Delta t}{St}\mathbf{U}^{(n)}.
\]
Next the Lagrangian force density is computed based on the new positions,
$\mathbf{X}^{(n)}$ as $\mathbf{F}^{(n)}=\mathcal{F}(\mathbf{X}^{(n)})$.
Finally, the Eulerian fields, $\mathbf{q}=\{\mathbf{f},\,\mu,\rho$\}
are computed using discrete$\delta$ function interpolation to the
Eulerian grid through equations of the form,
\[
\mathbf{q}^{(n)}=\sum_{\mathbf{X}^{(n)}\in\mathcal{L}}\mathbf{Q}^{(n)}(\mathbf{X}^{(n)})\,\hat{\delta}(\mathbf{x}_{h}-\mathbf{X}^{(n)},\,\omega).
\]

In the simulations we conduct, the primary direction of fluid flow
is in the $z$ direction. The height is governed by the $y$ coordinate
and width by the $x$ coordinate.

\begin{figure}
\begin{centering}
\subfloat[]{\begin{centering}
\includegraphics[width=0.3\paperwidth]{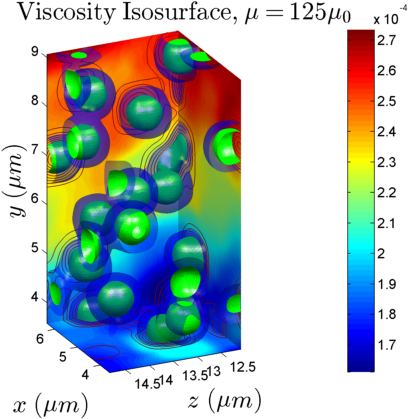}
\par\end{centering}

}\quad{}\subfloat[]{\begin{centering}
\includegraphics[width=0.3\paperwidth]{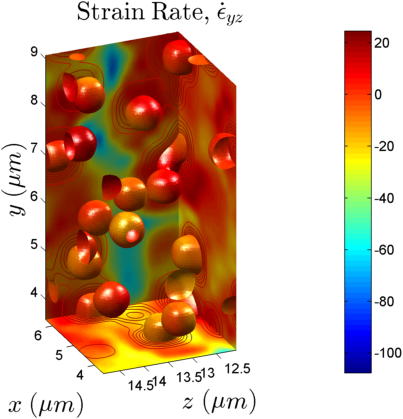}
\par\end{centering}

}
\par\end{centering}

\protect\caption{a) A viscosity isosurface is shown for a small section of a biofilm
used in simulation. The inner isosurface is $\mu=125\mu_{0}$ and
the outer transparent isosurface is at $\mu=50\mu_{0}$. Slices of
the $||\mathbf{u}||$ velocity field are shown as well. b) The $\dot{\epsilon}_{yz}$
component of the strain rate is plotted on the $\mu=125\mu_{0}$ viscosity
isosurface. Additionally, the strain rate and contours of viscosity
are shown in slice planes. For a single phase fluid, Newton's viscosity
law is $\sigma=\mu\dot{\epsilon}$. Although biofilms are not Newtonian
fluids, we still see that in areas of low viscosity, higher strain
rates are found and in areas of higher viscosity lower strain rates
occur. }
\end{figure}

\subsection{Numerical Verification and Convergence Properties}

In the first numerical verification result, we verify the accuracy
of the numerical projection method solver with no biofilm present
by comparing the numerical solution with an analytical solution. Since
there is no immersed structure, this is a test of the fluid solver
alone, and not the IBM method. For this test, the domain, $\Omega$
is chosen to be a rectangular solid that is periodic in the $x$ and
$z$ directions. From \cite{carlslaw1959conduction}, the following
boundary conditions for $y=0$ and $y_{L}$, 

\begin{eqnarray}
\frac{\partial P}{\partial y}=0\quad & \mathbf{u}|_{0}=\mathbf{0}\quad\mathbf{u}|_{y_{L}} & =\langle0,\,0,\,\sin\nu t\rangle
\end{eqnarray}
provide us with an analytic solution,

\begin{equation}
u_{z}(y,\, t)=\left|\frac{\sinh k\, y(1+i)}{\sinh k\, y_{L}(1+i)}\right|\sin\left(\nu\, t+\arg\left(\frac{\sinh k\, y(1+i)}{\sinh k\, y_{L}(1+i)}\right)\right)\,\,\,\,\, k=\left(\frac{\nu\rho}{2\mu}\right)^{1/2}.
\end{equation}
 The values of $P$, $u_{x}$ and $u_{y}$ are exactly zero in this
case. The values of $\rho$ and $\mu$ are set to $998\,\mathrm{kg}/m^{3}$and
$1\,\mathrm{Pa}\cdot s$ respectively and are homogenous across the
domain since no analytic solutions with variable density and viscosity
and the boundary conditions given above are known to the authors.
Convergence tests were conducted with frequencies $\nu=1\, Hz$ and
$\nu=100\, Hz$. In Table \ref{tab:ConvFactors} the absolute error,
temporal convergence factors, and spatial convergence factors are
listed. These quantities are calculated as described in Part I \cite{hammond2014variable}.

\begin{table}
\begin{centering}
\begin{tabular}{|c|c|c|c|}
\hline 
Frequency & Time & Space & Error $||\mathbf{u}^{h}-\mathbf{u}||_{\infty}$\tabularnewline
\hline 
$1\, Hz$ & 1.006 & 1.801 & $4.910\cdot10^{-10}$\tabularnewline
\hline 
$100\, Hz$ & 1.221 & 2.002 & $2.223\cdot10^{-5}$\tabularnewline
\hline 
\end{tabular}
\par\end{centering}

\protect\caption{\label{tab:ConvFactors}Spatial convergence tests were carried out
with grid spacings, $h$, set to 1/32, 1/64, and 1/128 and a time
step of $\Delta t=1/500/\nu$. Temporal convergence tests were done
with $\nu\Delta t$ set to 1/125, 1/250, and 1/500 and $\Delta x=1/64$.
Error is computed at $t=(0.2/\nu)\, s$. Convergence factors are computed
as $\rho(\Delta t)=\log\frac{||u(\Delta t/2)-u(\Delta t)||_{2}}{||u(\Delta t/4)-u(\Delta t/2)||_{2}}$
in time and by $\rho(h)=\log\frac{||u(h/2)-I_{h}^{h/2}u(h)||_{2}}{||u(h/4)-I_{h/2}^{h/4}u(h)/2||_{2}}$
where $I_{h}^{h/2}$ is an interpolation operator taking functions
from a grid with spacing $h$ to a grid with spacing $h/2$.}
\end{table}

Additionally, with the same boundary conditions as above, we tested
the convergence rates for simulations with a biofilm that possesses
variable density and viscosity. Temporal and spatial convergence factors
are shown in Table \ref{tab:Biofilm_Conv}. More detailed numerical
convergence results for this model with different boundary conditions
are shown in \cite{hammond2014variable}. In Table \ref{tab:Biofilm_Conv},
temporal convergence factors for the same fluid conditions and domain
as the analytical solution are listed. The pressure convergence rate
is not shown here since pressure variation only varies by about $\mathcal{O}(10^{-5})$
and thus is around the same order as the numerical errors observed
in the finite difference approximations used.

\begin{table}
\begin{centering}
\begin{tabular}{|c|c|c|c|c|}
\hline 
\multirow{2}{*}{Frequency} & \multicolumn{2}{c|}{Velocity, $||\mathbf{u}||$} & \multicolumn{2}{c|}{Position, $||\mathbf{X}||$}\tabularnewline
\cline{2-5} 
 & Time & Space & Time & Space\tabularnewline
\hline 
$49.91\, Hz$ & 0.983 & 1.105 & 1.022 & 0.952\tabularnewline
\hline 
$4.991\, Hz$ & 0.991 & 0.910 & 1.007 & 1.054\tabularnewline
\hline 
\end{tabular}
\par\end{centering}

\protect\caption{\label{tab:Biofilm_Conv}Convergence factors of hrIBM with biofilm.
For spatial convergence, $h$ was set to 1/32, 1/64, and 1/128 with
a time step of $\nu\Delta t=1/500$. To measure the temporal convergence
factors, $\nu\Delta t$ was set to 1/250, 1/500, and 1/1000. In both
cases, the boundary conditions described in Section \ref{sec:ConvergenceRate}
were used.}
\end{table}

\section{Experimental Validation Results\label{sec:Experimental-Validation}}

The material characterization of bacterial biofilms is a difficult
experimental task. It is usually not possible to grow biofilms large
enough for use in standard testing devices and, attempts to move a
biofilm from the environment it was grown in to a testing apparatus
may alter its structure \cite{pavlovsky2013insitu}. In Pavlovsky
et al. (2013) \cite{pavlovsky2013insitu}, a promising experimental
method of testing material properties of biofilms was developed. In
the experimental setup, a biofilm is grown in a parallel plate rheometer.
As the biofilm grows, it adheres to both the top and bottom plate
of the device. The top plate can then be rotated or repositioned vertically
and the stress and strain induced in the biofilm can be monitered.
These measurements can then be used to infer material properties of
the biofilm. Using the hrIBM model, we set up a simulation to reproduce
experiments described in Pavlovsky et al (2013). 

In order to reproduce the biofilm in simulation, 3D position data
sets obtained by high resolution microscopy of live biofilms are used
to initialize the positions of bacteria in the computational domain.
The experimental setup used to obtain these data sets are described
in Pavlovsky et al. (2015) \cite{pavlovsky2015effects} and Stewart
et al. \cite{stewart2013roleof}. Although the biofilm position data
sets that we use, which were obtained from the experiments described
in \cite{pavlovsky2015effects}, are not the ones grown and tested
in the bioreactor, they are from biofilms grown under similar physical
and nutrient availability conditions.  A key result seen from our
simulations is that the material properties computed by our model
of the different data sets are similar to each other. This indicates
that the material properties obtained through simulation must depend
on larger scale structural properties of the biofilm and may be treated
as bulk properties of the biofilm. For validation we compared bulk
properties measured by our model to experimental results. The methods
used to compute these quantities are discussed in the next subsections.

In Pavlovsky et al. (2013) \cite{pavlovsky2013insitu}, small amplitude
rheometry (SAR) is used to characterize the viscoelastic behavior
of \emph{S. epidermidis} biofilms. In SAR experiments, the upper plate
of the rheometer is rotated to induce a sinusoidal shear deformation
such that the average strain amplitude at the top of the biofilm is
a fixed value and the corresponding stress is measured. The strain
amplitude was set to $0.13$ at the outer radius of the rheometer
since this strain amplitude is found to be in a regime of primarily
linear and elastic mechanical behavior \cite{pavlovsky2013insitu}.
Using \ref{eq:DynamicModuli}, the dynamic moduli are computed at
a number of different frequencies of oscillation. 

Creep compliance testing is another characterization technique used
in Pavlovsky et al. (2013). In a creep compliance test, a constant
shear stress is applied to the biofilm through the top plate of the
rheometer. This induces a time dependent strain which can be measured. 

With the hrIBM model, we assume that for a small rectangular sample
of the biofilm that is not near the rotational center and, does not
border the outer boundary of the disc, the effects of cylindrical
geometry are negligible and the rotational motion can be approximated
as linear shear. This assumption greatly reduces the computational
expense of simulating the biofilm and simplifies the discretization
of the computational domain. This approximation is valid since the
stresses due to angular momentum are much less than those due to the
shearing motion of the plates. The size of the bioreactor used experimentally
is $40\, mm$ in diameter and approximately $250\,\mu m$ in height,
whereas, the computational domain is only $9-18\,\mu m$ in width
and length and $18-27\,\mu m$ in height. 

From Christensen \cite{christensen1982theory}, the shear stress and
strain, $\sigma_{\theta z}$ and $\epsilon_{\theta z}$, of an isotropic
viscoelastic cylinder undergoing small angle torsion are proportional
to $r$, the radial coordinate. Thus, if we choose to simulate some
subset of the bioreactor that is $30\,\mu m$ in width (this is larger
than in simulations we conduct) that is near, but not touching the
outer edge of the cylinder, at a radius of $15\, mm$ from the center,
the ratio of shear strain and shear stress exerted at the inner and
outer boundaries is approximately $(15\, mm)/(15\, mm+30\,\mu m)=0.998$.
Additionally, we see that in this case, the shear strain and shear
stress are not functions of the angular coordinate, thus approximating
the slightly curved domain as a rectangular solid should not alter
the physics of the problem. Of course, biofilms have far more complicated
material properties, however, at the scale of our simulations, this
result indicates that rectangular geometry and linear shear produces
an accurate approximation of the motion of the biofilm..

\subsection{Computation of Rheological Properties}

In order to compute the desired dynamic moduli, and compliance modulus
results, the stress and strain experienced by the biofilm during simulation
must be computed. The stress $\boldsymbol{\sigma}$, is decomposed
into a sum of stress due to the fluid motion, $\boldsymbol{\sigma}^{f}$,
and stress due to the straining of inter-bacteria connections within
the biofilm $\boldsymbol{\sigma}^{b}$. The total stress can then
be found as $\boldsymbol{\sigma}=\boldsymbol{\sigma}^{b}+\boldsymbol{\sigma}^{f}.$
Although each component of stress is computed separately during simulations,
distinct simulations cannot be used to individually test $\boldsymbol{\sigma}^{b}$
and $\boldsymbol{\sigma}^{f}$ since they are coupled. 

In order to calucate the strain $\boldsymbol{\epsilon},$ a set of
tracer particles is tracked throughout the simulation. Spatial derivatives
can then be calculated to obtain approximations of the strain. Additionally,
since only small amplitude strains are observed, the linear relation,
$\boldsymbol{\epsilon}=\frac{1}{2}(\nabla\mathbf{d}+\nabla^{T}\mathbf{d})$
is an accurate approximation of the strain for a displacement vector
$\mathbf{d}$. The derivatives needed to compute the strain are taken
with respect to the advected material coordinates.

Viscoelastic materials are often characterized through their time
dependent stress response to strain or their time dependent strain
response to stress. For a general viscoelastic material, given that
the stress and strain are sufficiently smooth functions of time, constitutive
relations between the stress and strain may be written in terms of
a convolution with viscoelasticity tensors as:

\begin{equation}
\sigma_{ij}(\mathbf{x},\, t)=\int_{-\infty}^{t}G_{ijkl}(\mathbf{x},\thinspace t-\tau)\,\frac{d}{d\tau}\epsilon_{kl}(\mathbf{x},\,\tau)\thinspace d\tau\label{eq:Relaxation_Const}
\end{equation}
\begin{equation}
\epsilon_{ij}(\mathbf{x},\, t)=\int_{-\infty}^{t}J_{ijkl}(\mathbf{x},\, t-\tau)\,\frac{d}{d\tau}\sigma_{kl}(\mathbf{x},\,\tau)\, d\tau\label{eq:Compliance_Const}
\end{equation}
where $\sigma_{ij}$ is the stress tensor, $\epsilon_{ij}$ is the
strain tensor, and $G_{ijkl}$ and $J_{ijkl}$ are fourth order viscoelasticity
tensors (see Christenson \cite{christensen1982theory} , \S 1 for
a derivation). In the literature, $\boldsymbol{G}$ is often called
the relaxation modulus and $\boldsymbol{J}$ is called the compliance
modulus. For linear, isotropic materials, the expression for $G{}_{ijkl}$
simplifies to $G_{ijkl}=\frac{1}{3}\left(G_{1}(t)-G_{2}(t)\right)\hat{\delta}_{ij}\hat{\delta}_{kl}+\frac{1}{2}G_{1}\left(\hat{\delta}_{ik}\hat{\delta}_{jl}+\hat{\delta}_{il}\hat{\delta}_{jk}\right)$,
where $\hat{\delta}_{mn}$ is the Kronecker delta function and Einstein
summation notation is used. The two functions, $G_{1}(t)$ and $G_{2}(t)$
correspond respectively to shear and dilatational stresses. Analogous
expressions exist for the compliance tensor. Although the viscoelastic
moduli are spatially heterogeneous, we believe that more meaningful
results are obtained in the mean field, or spatially averaged, time
dependent values for $\boldsymbol{\epsilon}$, $\boldsymbol{\sigma}$,
$\boldsymbol{G}$, and $\boldsymbol{J}$. These quantities depend
less on the exact configuration of bacteria in a biofilm and behave
more like bulk material parameters that can be measured experimentally.
Although the interconnected links used to model the connections between
adjacent bacteria each individually introduce anistropy into the model,
under the conditions of our simulations, the overall behavior of the
biofilm is not highly anisotropic.

\subsubsection{Computation of Strain}

Although a single phase Newtonian fluid will behave viscously (i.e.,
the stress only depends on the strain rate, not strain itself), in
a biofilm the fluid component is influenced by the elastic components
of the biofilm and thus the stress state in a biofilm depends directly
on the strain (along with the strain rate). In order to compute the
strain, the displacement field must be computed. The displacement,
$\mathbf{d}$ of a particle located at $\mathbf{x}_{0}$ at time,
$t=0$ in a material undergoing deformation can be found by solving
the following ODE:
\begin{equation}
\frac{\partial}{\partial t}\mathbf{d}(\mathbf{x}_{0},\, t)=\int_{\Omega}\mathbf{u}(\mathbf{x},\, t)\delta(\mathbf{x}-\mathbf{d}(\mathbf{x}_{0},\, t)-\mathbf{x}_{0})\, d\mathbf{x};\quad\mathbf{d}(\mathbf{x}_{0},\,0)=\mathbf{0}.
\end{equation}
In the biofilm simulations, ``tracer'' particles with positions
denoted by $\mathbf{S}(x,\, y,\, z)$, are initialized at heights
$y_{L}-\gamma$, $y_{L}-\gamma-h$, and $y_{L}-\gamma-2h$, near the
top of the biofilm at $t=0.$ At each time step, the positions of
the tracers are updated using the same $\delta$ function interpolation
used to update the bacteria positions. With these tracers, the deformation
of the biofilm can be tracked throughout the simulation.

In the simulations, the $\epsilon_{yz}$ component of strain is needed
at the upper boundary of the domain. Therefore, the tracers are initialized
near the top of the domain in three vertically aligned layers. This
is done to make the numerical approximation of derivatives of the
form $\partial d_{z}/\partial y$ easier . With the initial arrangement
of tracters in vertically aligned layers, the centered finite difference
approximation
\begin{equation}
\epsilon_{yz}(\mathbf{S},\, t)\approx\frac{1}{2}\left(\frac{\partial d_{z}}{\partial y}+\frac{\partial d_{y}}{\partial z}\right)\approx\frac{1}{2}\left(\frac{\partial d_{z}}{\partial y}\right)\approx\frac{1}{2}\left(\frac{1}{2}\frac{d_{z}(\mathbf{S}(y)-d_{z}(\mathbf{S}(y-h))}{S_{y}(y)-S_{y}(y-h)}+\frac{1}{2}\frac{d_{z}(\mathbf{S}(y+h)-d_{z}(\mathbf{S}(y))}{S_{y}(y+h)-S_{y}(y)}\right)
\end{equation}
can be used to approximate the strain. The reported value of $\epsilon_{yz}$
at each time step is then the average of the strains calculated over
each tuple of tracers. Since the entire upper plate moves at a single
velocity at any given time, $\partial d_{y}/\partial z$ is negligible
in this case, whereas in general, this term is required to compute
the shear strain.

\subsubsection{Computation of Stress Induced by Fluid Motion}

From Newton's viscosity law the $\sigma_{yz}^{f}$ component of stress
can be found as 
\begin{equation}
\sigma_{yz}^{f}=\mu(\mathbf{x})\left(\frac{\partial u_{z}}{\partial y}+\frac{\partial u_{y}}{\partial z}\right).
\end{equation}
Since the velocity field is already known from solving the Navier-Stokes
equations at each time step, the relevant derivatives can be approximated
by finite difference approximations. As with the strain calculation,
the second term, $\partial u_{y}/\partial z$, is zero since the $y$
velocity on the entire top plate of the rheometer is zero. The reported
value of $\sigma_{yz}^{f}$ at each time step is then found by spatially
averaging over the top $2.5\,\mu m$ of the domain. This is done instead
of just averaging over the very top of the domain in case there are
numerical boundary layers in the fluid flow field near the boundary.
Boundary layers of thickness $\mathcal{O}(\sqrt{\mu\Delta t/\rho})$
are known to sometimes arise in projection method based fluid solvers
\cite{brown2001accurate,orszag1986boundary}.To mitigate this problem,
we use boundary conditions that do not cause this issue in the constant
density and viscosity case.

\subsubsection{Computation of Stress Induced by the Biofilm Configuration}

In order to compute the force exerted by the biofilm connections on
the top plate, we integrate the Eulerian force field induced by bacteria
adhered to the top plate. To determine if a bacteria is adhered, we
choose a distance, $\gamma=0.4\,\mu m$ from the top plate, and assume
that each bacteria with $y$ coordinate in the interval $[y_{L}-\gamma,\, y_{L}]$
is adhered to the top, and that its $z$-component of velocity is
fixed to be that of the upper plate. For these bacteria, any force
applied on them by spring-like connections to other bacteria behaves
like a force exerted by the biofilm on the upper plate instead of
on the bulk fluid. The sum of these forces is used to compute the
stress induced by the spring-like connections by means of Cauchy's
traction law,
\begin{equation}
\boldsymbol{\sigma}^{b}\mathbf{n}=\frac{\mathbf{F}^{b}}{A}.
\end{equation}
The outward unit normal, $\mathbf{n}$, is $(0,\,1,\,0)$ in this
case since the top plate is parallel to the $xz$ plane. The force,
$\mathbf{F}^{b}$ is found by integrating the Eulerian force density
field that would be generated by the biofilm nodes adhered to the
top plate. Additionally, since we are interested in the applied shear
stress, $\sigma_{zy}^{b}$, this can be found as -$\frac{F_{z}^{b}}{A}$.
Note that $\gamma$ was chosen arbitrarily, however we observed that
with $\gamma=0.7\,\mu m$ the results were not significantly different.

\subsection{Shear Moduli, $G^{\prime}$and $G^{\prime\prime}$\label{sub:Shear-Moduli}}

When a nearly isotropic material is subjected to an oscillatory displacement
field with frequency $\nu$, we may write the strain as $\epsilon(t)=i\nu\thinspace\epsilon_{0}e^{i\nu t}$,
where $i$ is the imaginary unit and $\epsilon_{0}$ is the strain
amplitude. For cases where the strain is primarily only shear strain
equation (\ref{eq:Relaxation_Const}) gives, $\boldsymbol{\sigma}(\nu)\approx G_{1}^{\ast}(\nu)\boldsymbol{\epsilon}(\nu)$
where $G_{1}^{\ast}(\nu$) is related to the Fourier transform in
time of $G_{1}(t)$. In general $G_{1}^{\ast}(\nu)$ is a complex
valued function. Breaking the complex shear modulus into its real
and imaginary components, $G_{1}^{\ast}(\nu)=G^{\prime}(\nu)+iG^{\prime\prime}(\nu)$;
and given a strain amplitude $\epsilon_{0}(\nu)$ and stress amplitude
$\sigma_{0}(\nu)$ (in Pascals), 
\begin{equation}
G^{\prime}(\nu)=\frac{\sigma_{0}(\nu)}{\epsilon_{0}(\nu)}\cos\delta(\nu),\qquad G^{\prime\prime}(\nu)=\frac{\sigma_{0}(\nu)}{\epsilon_{0}(\nu)}\sin\delta(\nu).\label{eq:DynamicModuli}
\end{equation}
Here, $\delta(\nu)$ is known as the loss angle, measured in radians
at frequency $\nu$. In the literature, $G^{\prime}(\nu)$ and $G^{\prime\prime}(\nu)$
are often referred to as the storage and loss moduli. They correspond
to the elastic and viscous components of a viscoelastic stress strain
relationship. 

\begin{figure}
\begin{centering}
\includegraphics[width=0.4\columnwidth]{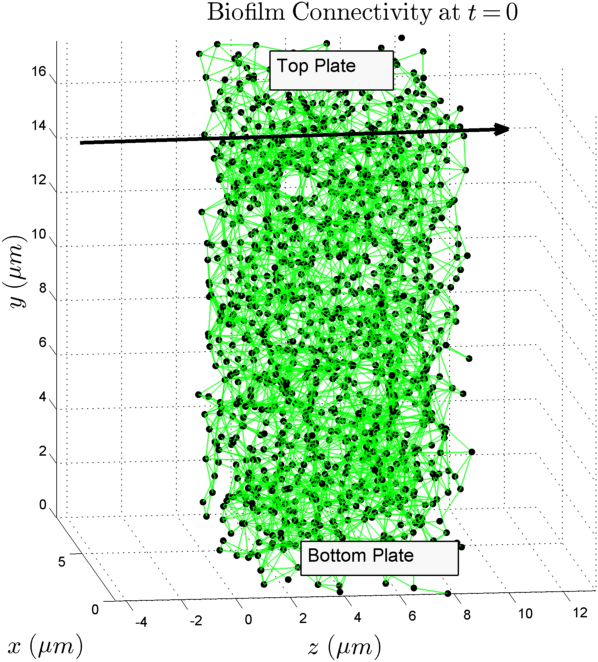}\enskip{}\includegraphics[width=0.4\columnwidth]{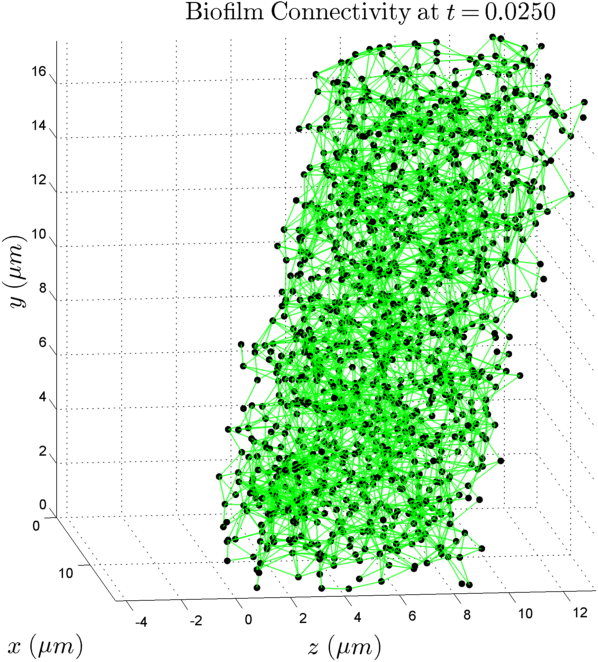}\\
\includegraphics[width=0.4\columnwidth]{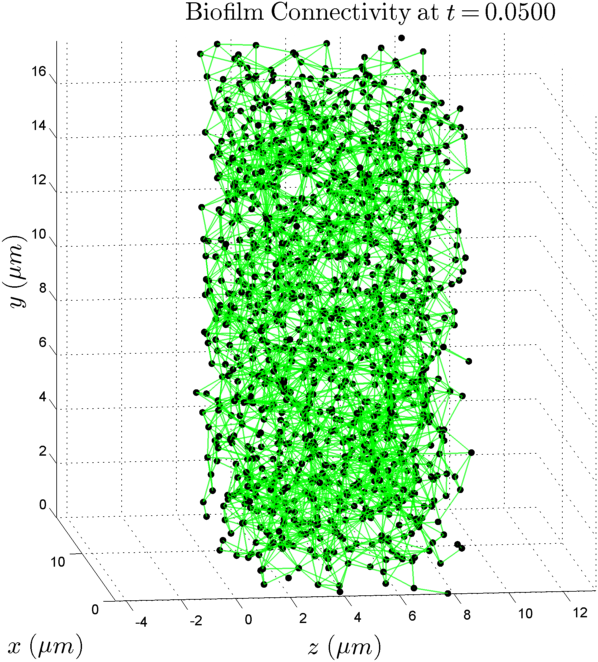}\enskip{}\includegraphics[width=0.4\columnwidth]{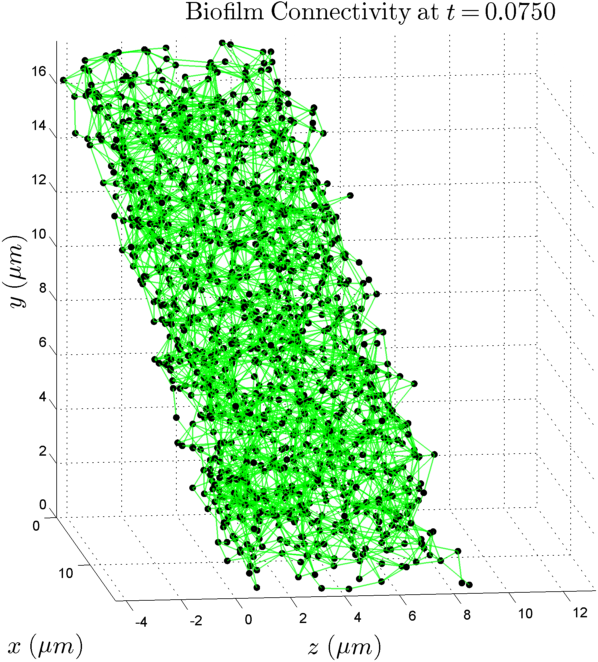}
\par\end{centering}

\protect\caption{\label{fig:SAR_Deformation}Starting at $t=0$ on the left, the images
show how the biofilm is moved as the top plate oscillates. Dots are
bacteria locations and lines indicate viscoelastic connections. In
the simulations the domain is periodic in the $x$ and $z$ directions.
The periodicity is not shown here since it makes it more difficult
to visualize the effect of deformation on the biofilm.}
\end{figure}

Taking the domain to be a rectangular solid, oriented as shown in
Figure \ref{fig:Biofilm_Geom}, we assume that all fields are periodic
in the $x$ and $z$ directions. We use the following boundary conditions:
\begin{equation}
\begin{array}{ccc}
\left.\frac{\partial P}{\partial y}\right|_{y=0,y_{L}}=0\quad & \mathbf{u}(x,0,z,t)=0\quad & \mathbf{u}(x,y_{L},z,t)=(0\:,0\:,u_{b}(t)\,).\end{array}
\end{equation}
Along the top boundary, we set the $z$ velocity to be 

\begin{equation}
u_{b}(t)=\epsilon_{0}\frac{(e^{2\nu t}-1)((e^{4\nu t}-1)\cos\nu t+8e^{2\nu t}\sin\nu t)}{(1+e^{2\nu t})^{3}}.
\end{equation}
This particular function is chosen since it is continuous, at $t=0$,
$u_{z}=0$, and because it converges to within 0.001 of $\epsilon_{0}\cos\nu t$
within half an oscillation, reducing the amount of time needed to
run simulations. To initialize the bacteria positions, we take a $9\,\mu m\times27\,\mu m\times9\,\mu m$
subset of a $30\,\mu m\times30\,\mu m\times10\,\mu m$ bacteria position
data field obtained experimentally. This data is also used in the
initialization of the viscosity and density fields present in the
biofilm. We believe that setting the internal forces to zero at the
start is reasonable since experimental results from SAR under both
compression and tension yielded similar results. In Figure \ref{fig:SAR_Deformation},
the deformation induced by an oscillatory shearing motion is depicted.

In order to tune our model to the experimental data, we adjusted the
spring constant, $k_{ij}$ used in Hooke's Law and the distance by
which we allow any two biofilm nodes to be connected by at time $t=0$.
For a spring connecting two points in space, Hooke's Law can be written
as 
\begin{equation}
\mathbf{F}_{ij}=k_{ij}\,\Lambda_{ij}(\mathbf{X},\,\mathbf{X}_{0})\,(\mathbf{X}_{i}-\mathbf{X}_{j}),
\end{equation}
with
\begin{equation}
\Lambda_{ij}(\mathbf{X},\,\mathbf{X}_{0})=\frac{||\mathbf{X}_{i}(t)-\mathbf{X}_{j}(t)||-||\mathbf{X}_{i}(0)-\mathbf{X}_{j}(0)||}{||\mathbf{X}_{i}(t)-\mathbf{X}_{j}(t)||}.
\end{equation}
Following Hammond et al. \cite{hammond2014variable}, we choose each
$k_{ij}$ to be a force constant, $F_{max}$ divided by the initial
separation of bacteria $i$ and $j$. Since the immersed boundary
method requires a Lagrangian force density, we then divide $\mathbf{F}_{ij}$
by the Lagrangian volume element, $d_{0}^{3}$. Additionally, we note
that in Dan Vo et al. \cite{danvo2010anexperimentally} and Peskin
\cite{peskin2002theimmersed}, an identical constitutive relation
is derived from the starting point of energy functionals in which
the force density is found by taking a Fr\'echet derivative of an
energy functional.

\begin{figure}
\begin{centering}
\includegraphics[width=0.5\columnwidth]{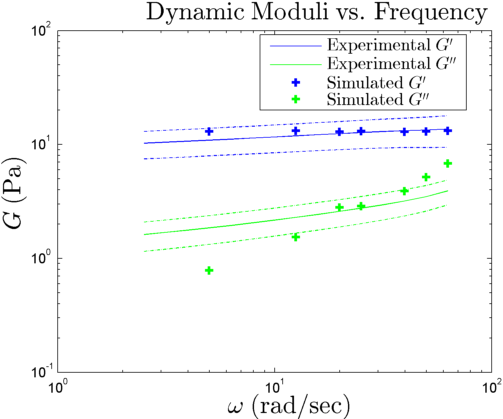}
\par\end{centering}

\protect\caption{\label{fig:G1_G2_vs_w} A comparison between experimentally measured
results for $G^{\prime}$ and $G^{\prime\prime}$ is shown in comparison
the simulation results. For these results, the force constant was
$F_{max}=1.3223\cdot10^{-9}$ , the connection distance between bacteria
was $1.62\mu m$, no damping was used, $\mu_{b}=250\mu_{0}$, and
$\rho_{b}=1.12\rho_{0}$. The dashed lines indicate the experimental
error range.}
\end{figure}

In Figure \ref{fig:G1_G2_vs_w} we depict the frequency dependence
of $G^{\prime}$ and $G^{\prime\prime}$ . From these results, it
is clear that our model fits experimental data on $G^{\prime}$ quite
well. For $G^{\prime\prime}$ the fit is not as strong, although we
still do see that many of the results from simulation are within the
range of experimental error. We observe that in fact the slope of
$G^{\prime\prime}$ is steeper than the experimental measurements.
Although at this time, the cause of this difference is unknown, it
is possible that extensions such as those discussed in Section \ref{sec:Future-Directions}
may correct this. 

\subsection{Creep Compliance Measurements $J(t)$}

Creep compliance is a measure of how a material deforms over time
in response to an applied stress. Experimentally, the compliance is
measured by using the rheometer to apply a step change to the shear
stress on the upper plate and observing the resultant shear strain.
A step change in stress can be written as, $\bar{\sigma}(t)=\sigma_{0}H(t)$
where $\sigma_{0}$ is the magnitude of the step change and, $H(t)$
is the Heaviside step function. For a linear isotropic material with
$\sigma_{yz}(t)=\bar{\sigma}(t)$, the integral in Equation (\ref{eq:Compliance_Const})
simplifies to
\begin{equation}
\epsilon_{yz}(t)=\sigma_{0}J_{1}(t).
\end{equation}

From a physical standpoint, most of the bacteria and the bulk of the
fluid are only effected by the step change in stress after the stress
propogates vertically through the biofilm. However, in the portion
of the biofilm adjacent to the upper plate the effect of a change
in stress is instantaneous. Thus, we can write a force balance between
the forces in the biofilm, the acceleration of the top plate, and
the applied force on the top plate. This leads to an impulse boundary
condition which specifies the velocity at the top plate  The boundary
condition can be written as
\begin{equation}
\left.\frac{d}{dt}u_{z}\right|_{y=H}=(\rho V)^{-1}(\sigma_{0}-\sigma^{b}-\sigma^{f})A
\end{equation}
 where $\sigma^{b}$ and $\sigma^{f}$ are the stress exerted by the
fluid and the springs in the biofilm at the top plate, $A$ is the
area of the upper plate, and $(\rho V)$ is the mass associated with
the top plate of the rheometer. We assume that this mass is equivalent
to the mass of the top $2.4\,\mu m$ of the biofilm where bacteria
are adhered to the top plate.

\begin{figure}
\begin{centering}
\includegraphics[width=0.5\columnwidth]{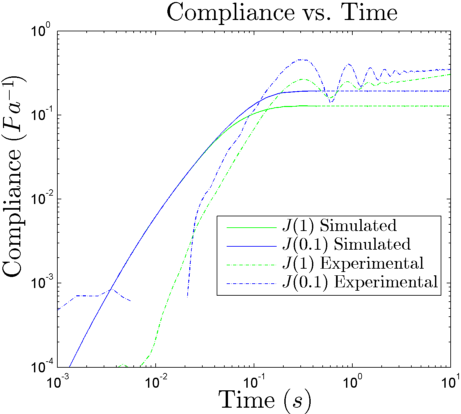}
\par\end{centering}

\protect\caption{\label{fig:Compliance}The time dependence of $J(t,\,\sigma_{0})$
is shown above for $\sigma_{0}=0.1\, Pa$ and $\sigma_{0}=1\, Pa$.
It can be seen that the compliance levels out at a similar value as
the experimental result, but levels out faster than in the experiments.
For the compliance simulations, we used $F_{max}=2.9091\cdot10^{-9}$
instead of $F_{max}=1.3223\cdot10^{-9}$ due to some numerical stability
issues.}
\end{figure}

Numerically, this boundary condition can be written as
\begin{equation}
u_{z}^{(n+1)}|_{y=H}=u_{z}^{(n)}|_{y=H}+\left(\frac{\Delta t}{\rho_{0}u_{0}L\sum_{\mathcal{G}_{h}}\rho_{ijk}}\right)(\sigma_{0}-\sigma^{b}-\sigma^{f}).\label{eq:Compliance_BC}
\end{equation}
In (\ref{eq:Compliance_BC}), $\rho_{0}$ and $\mu_{0}$ are the density
and viscosity of water, and $L$ is the characteristic length (in
this case $10\,\mu m$). In numerical experiments, rather than immediately
impose a step in stress at time 0, we add a mollifier, $\left(\frac{2}{1-e^{-\alpha t}}-1\right)$
on the applied stress, $\sigma_{0}$ to mitigate any possible numerical
instabilities associated with a discontinuous boundary condition.
In this case, $\alpha=200$ is chosen to be large so that the applied
stress approaches its equilibrium value within 0.1 seconds. This is
reasonable because very short time compliance behavior is not generally
experimentally measurable, and also a step in the stress may not actually
occur instantaneously from the perspective of a very short time scale.
Results from simulations using this boundary condition andtwo values
of $\sigma_{0}$ are shown in Figure \ref{fig:Compliance}.

Presently, we do not propose mechanisms by which the spring-like connections
may break and reconnect, although it is possible to incorporate such
a model into our current framework as discussed in \cite{hammond2014variable}.
Thus, the long term behavior of $J_{1}(t,\,\sigma_{0})$ which likely
depends on the gradual redistribution of the connectivity of the biofilm
is not expected to be captured. Therefore, we do not simulate beyond
0.1 seconds. In the experimental results, there are some damped oscillations
present in the creep compliance experiments. In Pavlovsky et al. (2013),
these oscillations are found to be related to intertial effects from
the rheometer itself and thus are not expected or observed in our
simulations.

\subsection{Similarity of Material Properties Between Different Bacteria Position
Data Sets}

\begin{figure}
\begin{centering}
\subfloat[]{\begin{centering}
\includegraphics[scale=0.4]{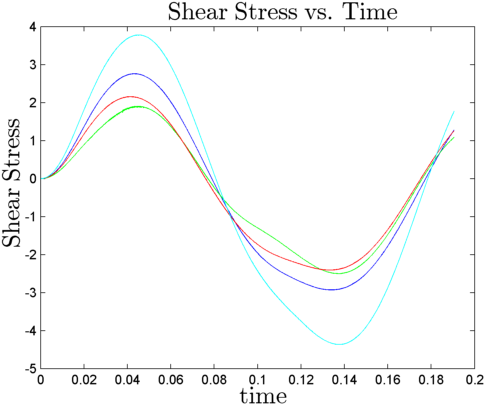}
\par\end{centering}

}\subfloat[]{\begin{centering}
\includegraphics[scale=0.4]{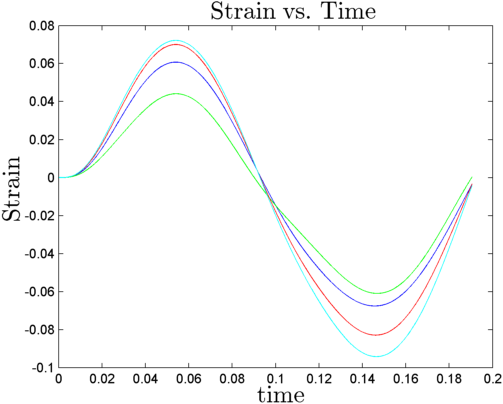}
\par\end{centering}

}
\par\end{centering}

\protect\caption{\label{fig:Stress_Strain_Similarity}These graphs show the stress
vs. time and strain vs. time for 4 different biofilm samples. The
samples are all $18\mu m\times27\mu m\times9\mu m$ size and contain
approximately 2000 bacteria positions. }
\end{figure}

In Dzul et al. \cite{dzul2011contribution}, the spatial statistics
of bacteria in a biofilm are studied. We show here that from data
sets that have similar spatial distributions of nearest neighbor connections
between bacteria, similar bulk property measurements are obtained.
In our simulations, we take blocks of experimental data that are $18\,\mu m$
wide, $9\,\mu m$ long, and $27\,\mu m$ high and compute the dynamic
moduli of each block at a fixed frequency. . The graphs in Figure
\ref{fig:Stress_Strain_Similarity} show the stress and strain of
four different biofilms over one period of oscillation. 

From figure \ref{fig:Stress_Strain_Similarity} and Table \ref{tab:G1G2Results},
we see that in three of the four biofilm data sets, similar results
are obtained.. We also note that the mean standard error in the experimental
results for this particular test was 3.9791 for $G^{\prime}$ and
0.8836 for $G^{\prime\prime}$. We also suspect that if larger data
sets were used, even better agreement would be seen in the computed
values of $G^{\prime}$ and $G^{\prime\prime}$. 

The importance of this section is in verifying that the properties
we are validating can be considered as bulk properties. Since three
out of the four data sets provided results within the experimental
error deviation, we believe this to be strong evidence the properties
we measure are bulk properties.

\begin{table}
\begin{centering}
\begin{tabular}{|c|c|c|c|}
\hline 
Simulation & 1 & 2 & 3\tabularnewline
\hline 
$G^{\prime}(\nu=49.91)$ & 13.06 & 9.18 & 10.03\tabularnewline
\hline 
$G^{\prime\prime}(\nu=49.91)$ & -5.16 & -4.44 & -3.92\tabularnewline
\hline 
$\delta(\nu=49.91)$ & 0.376 & 0.451 & 0.373\tabularnewline
\hline 
\end{tabular}
\par\end{centering}

\protect\caption{\label{tab:G1G2Results} Results for $G^{\prime}(\nu)$, $G^{\prime\prime}(\nu)$
and $\delta(\nu)$ are shown at two frequencies for 3 different biofilm
coordinate data sets with $\nu=49.91\, rad/s$. These results show
that the physical properties measured here do not depend solely on
the exact microstructure of the biofilm, but on sometype of more large
scale organization of the bacteria positions in space.}
\end{table}

\subsection{In-Stream Tumbling of Biofilm Fragment}

Bacterial structures exhibit a diverse range of interaction with fluid
flow. One such interaction is the tumbling motion of aggregates in
shear flow. To simulate this effect, we conduct simulations in which
there are no bacteria anchored along the plates of the domain. Instead,
an aggregate of bacteria is located near the middle of the computational
domain and the upper and lower plates move in opposite directions
as 
\begin{equation}
u_{z}(x,\, y_{L},\, z)=10^{-3}\left(\frac{1}{1+e^{-t}}-\frac{1}{2}\right)=-u_{z}(x,\,0,\, z).
\end{equation}
The boundary velocities are scaled by $10^{-3}$ so as to limit the
shear stress so that the aggregate is not simply torn apart. To ensure
that the biofilm is not attached to the plates and is sufficiently
far from the plate to induce a rotating, or tumbling motion, these
simulations only include bacteria that are greater than $8.8\,\mu m$
from either plate at the start. With the physical parameters we use,
the bacteria aggregation rotates and is deformed by the fluid shear
forces exerted by the fluid \cite{byrne_postfragmentation_2011}.
\begin{figure}
\begin{centering}
\includegraphics[width=0.32\textwidth]{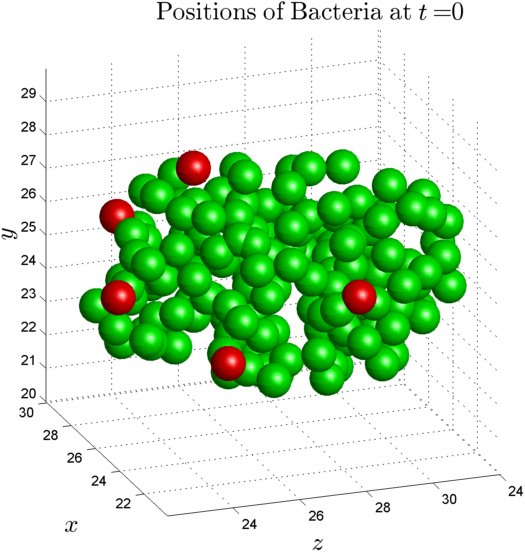}\quad{}\includegraphics[width=0.3\textwidth]{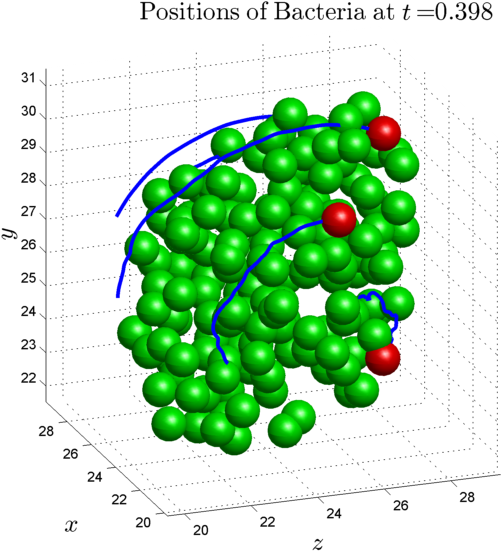}\quad{}\includegraphics[width=0.3\textwidth]{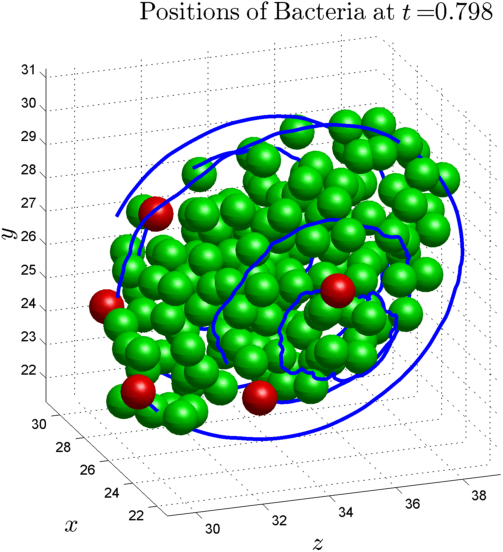}
\par\end{centering}

\protect\caption{\label{fig:Biofilm-aggregate-suspended}Biofilm aggregate suspended
in shear flow rotate over time. Snap shots shown at 0, 0.4 and 0.8
seconds into the simulation. The flattening of the ellipse in response
to the shear flow can be distinguished between the first and third
figure. Several bacteria are marked red to help show the rotation
of the aggregate. The blue lines indicates the trajectories of the
marked cells relative to the center of mass of the aggregate. Distances
are in 10s of micrometers.}
\end{figure}

In Blaser et al. \cite{blaser2002forces}, analytical results on the
frequency at which a solid ellipsoid will rotate in shear flow are
provided. For an ellipsoid with axis aligned with the direction of
fluid motion, the rotational frequency is found as 
\[
T=\frac{2\pi(a_{1}^{2}+a_{2}^{2})}{a_{1}a_{2}\tau}
\]
 where $\tau$ is the shear rate and $a_{1}$and $a_{2}$ are the
principle axes of the ellipse undergoing rotation. In our simulation
(see Figure \ref{fig:Biofilm-aggregate-suspended}), we show that
a bacterial aggregate approximated as a hydrodynamically equivalent
ellipse will rotate at a frequency similar to the theoretically expected
result. The rotational frequency of the aggregate is found by computing
the average frequency of rotation of bacteria in the $yz$ plane about
the center of mass of the aggregate. In the simulations, we observed
a frequency of approximately 0.76 seconds for an aggregate approximated
by an ellipse with major axis $a_{1}=0.167\,\mu m$ and first semimajor
axis $a_{2}=0.136$, and with shear rate of $20\, s^{-1}$ containing
54 bacteria. The theoretical result in this case is 0.641 seconds.
These results are shown in Table \ref{tab:ComparisonRotation}. Although
in our case these results are not as concrete of a metric as the dynamic
modulus and compliance results, rotational frequency has been used
for model validation in fields such as red blood cell modeling (see
\cite{fai2013immersed}).

\begin{table}
\begin{centering}
\begin{tabular}{|c|c|c|c|c|}
\hline 
Major axis, $a_{1}$ & First minor axis, $a_{2}$ & Theoretical Period & Observed Period & Relative Error\tabularnewline
\hline 
\hline 
$1.67\mu m$ & $1.36\mu m$ & $0.64s$ & $0.76s$ & $+18.8\%$\tabularnewline
\hline 
$4.94\mu m$ & $2.74\mu m$ & $0.59s$ & $0.67s$ & $+13.5\%$\tabularnewline
\hline 
\end{tabular}
\par\end{centering}

\protect\caption{\label{tab:ComparisonRotation}Comparison of theoretical and simulated
rotational results}
\end{table}

\section{Discussion and Future Directions\label{sec:Future-Directions}}

Although we see that the hrIBM model provides a versatile means of
simulating biofilms and can accurately capture some of the experimentally
observed behavior of biofilms, there is still room to extend the model
to allow for more general modeling of biofilms. An important area
of research in biofilm studies is developing an understanding of the
fracture mechanics and dynamics of biofilms. One way to model fracture
dynamics in our biofilm model would be the inclusion of a stochastic
model governing the connectivity of the viscoelastic links in the
biofilm. Thus, we could define a probability based on the stress and
strain in each link and the proximity of each pair of connected bacteria
to allow for reconfiguration of the connectivity of the biofilm over
time. It is also possible to model the viscoelastic properties of
the biofilm by adjusting the constitutive model that is used to provide
the Lagrangian force based on the nodal configuration of the bacteria.
Modeling the changing connectivity of a biofilm was explored in \cite{alpkvist2006threedimensional,bottino1998modeling,sudarsan2015simulating}.

Another area that could be explored is the shape of the discrete $\delta$
function used to approximate each bacteria and its associated viscosity
halo. It is possible that adjusting this function may allow for more
accurate modeling of the mass displacement induced by the bacteria
bodies in the bulk fluid. Adjustments to the $\delta$ function may
also allow for the inclusion of nonspherical bacteria into the model.
It would be interesting to see if similar results are obtained for
bacteria that are different shapes. 

One possible difficulty in adjusting the discrete $\delta$ function
is the preservation of mass in the model. In Equation (\ref{eq:Continuity_Eq}),
there is no density dependence as would normally be seen with the
Navier-Stokes equations for a variable density system. For the original
IBM, this is in fact exact as described in \cite{peskin2002theimmersed}.
For the hrIBM, there is an error however, error term is expected to
be small in our situation since $\rho(\mathbf{x},\, t$) only varies
by $12\%$ over the domain (density of bacteria is not highly variable),
all simulations are at low Reynold's numbers, and the high viscosity
gradients which overlap where density gradients occur make it difficult
for fluid to rapidly travel down a density gradient. Although we believe
that this approximation is accurate in the simulations we conduct,
this error term may increase in situations with higher Reynold's numbers.
This is an area that may be further developed in future works.

Another area of potential improvement is the numerical methods. Currently
the scheme is $\mathcal{O}(\Delta t$) and $\mathcal{O}(h$). In future
work, the Crank-Nicholson time-stepping scheme may be used since,
at least in the constant viscosity and density case, it could lead
to $\mathcal{O}(\Delta t^{2})$ convergence as shown in Brown et al.
\cite{brown2001accurate}. Another approach is to use a predictor-corrector
type of method such as those described in \cite{fai2013immersed}
and \cite{peskin2002theimmersed}. However, even without heterogeneous
material properties, obtaining $\mathcal{O}(\Delta t^{2})$ convergence
in the overall IBM is more complicated and also depends on properties
of the discrete Dirac $\delta$ function. A detailed discussion can
be found in Liu and Mori \cite{liu2012properties,liu2014lpconvergence}.
The development of efficient $\mathcal{O}(\Delta t^{2})$ IBM schemes
with heterogeneous material properties is still an area of active
research.

\section{Conclusions}

Based on the experimental results shown above and in \cite{hammond2014variable},
we show that our biofilm model, Equations (\ref{eq:N-S})-(\ref{eq:visc_transfer}),
can be used to determine bulk material properties of bacterial biofilms.
In particular, we show that the model yields close agreement with
experimental results from \cite{pavlovsky2013insitu} in which the
bacterium \emph{S. epidermidis} was grown in a bioreactor and characterized
using a parallel plate rheometer. We also show that suspended aggregates
of bacteria in shear flow rotate with a similar period as a hydrodynamically
equivalent ellipse. Another development is the uniformity of bulk
material properties over different experimental data sets that possess
similar spatial statistics. An important step in obtaining these results
was the computation of bulk material properties from simulations.
To our knowledge, our model is the first that can compute bulk material
properties of biofilms based on direct simulation of both microscale
connectivity of the biofilm and the heterogeneous rheology of the
ECM. 

We also acknowledge a number of new research directions and extensions
that can be done to improve results and also to allow for more flexible
modeling of biofilms in different scenarios than what we have considered
here.

\section{Acknowledgements}

This work was supported by the National Science Foundation grants
PHY-0940991 and DMS-1225878 to DMB, and PHY-0941227 to JGY and MJS,
and by the Department of Energy through the Computational Science
Graduate Fellowship program, DE-FG02-97ER25308, to JAS. This work
utilized the Janus supercomputer, which is supported by the National
Science Foundation (award number CNS-0821794), the University of Colorado
Boulder, the University of Colorado Denver, and the National Center
for Atmospheric Research. The Janus supercomputer is operated by the
University of Colorado Boulder. 

\bibliographystyle{siam}
\bibliography{Biofilm_Validation_Paper}

\end{document}